\documentclass[10pt]{article}

\usepackage {amsfonts}
\usepackage {amssymb}
\usepackage {amsmath}
\usepackage {hyperref}
\usepackage{graphicx}
\usepackage{amsthm}
\usepackage {color}
\usepackage{calc}

\makeatletter
\newlength{\temp@wc@width}
\newlength{\temp@wc@height}
\newcommand{\widecheck}[1]{%
  \setlength{\temp@wc@width}{\widthof{$#1$}}%
  \setlength{\temp@wc@height}{\heightof{$#1$}}%
  #1\hspace{-\temp@wc@width}%
  \raisebox{\temp@wc@height+2pt}[\heightof{$\widehat{#1}$}]%
     {\rotatebox[origin=c]{180}{\vbox to 0pt{\hbox{$\widehat{\hphantom{#1}}$}}}}%
}
\makeatother

\newcommand{\A}{\forall}

\newcommand{\Ext}{\mathbf{E}}
\newcommand{\Sh}{\mathbf{S}}
\newcommand{\Lr}{\mathbf{L}}

\newtheorem{Ps}{Proposition}[section]

\newtheorem{T}{Theorem}[section]
\newtheorem{Lm}{Lemma}[section]
\newtheorem{D}{Definition}[section]

\newtheorem{C}{Corollary}[section]

\title{Concave majorant of stochastic processes and  Burgers turbulence}
\author{Rapha\"el \textsc{LACHI\`EZE-REY}\thanks{Equipe Probabilit\'es et statistiques, UFR Math\'ematiques, Univ. Lille 1, raphael.lachieze-rey@math.univ-lille1.fr}}

\date{}
\begin{document}
\maketitle
\vspace{2cm}

\textbf{Abstract}
The asymptotic solution of the inviscid Burgers equations with initial  potential $\psi$ is closely related to the convex hull of the graph of $\psi$.

In this paper, we study  this convex hull, and more precisely its extremal points, if $\psi$ is a stochastic process. The times where those extremal points are reached, called extremal times, form a negligible set for L\'evy processes, their integrated processes, and It\^o processes. We examine more closely the case of a L\'evy process with bounded variation. Its extremal points are almost surely countable, with accumulation only around the extremal values. These results are derived from the general study of the extremal times of $\psi+f$, where $\psi$ is a L\'evy process and $f$ a smooth deterministic drift.

 These results allow us to show that, for an inviscid Burgers turbulence with a compactly supported initial potential $\psi$, the only point capable of being   Lagrangian regular   is the time $T$ where $\psi$ reaches its maximum, and that is indeed a regular point iff $0$ is regular for both half-lines. As a consequence, if the turbulence occurs on a non-compact interval, there are a.s. no Lagrangian regular points.

\vspace{1cm}


\textbf{AMS}: 60G17, 60D05, 35Q35, 35R60, 60J99.\\

\section*{Introduction and notation}

There are multiple reasons to study the convex hull $C_{X}$ of a random process $X$. From a geometric point of view, it is a good indicator of the graph  structure, both locally and globally, and captures some features of the regularity of $X$. Groenboom \cite{G} has first studied this topic if $X$ is a standard brownian motion and $I=\mathbb{R}^+$, establishing the essential results, and Pitman \cite{P} gave a simple and exhaustive description of the concave majorant, which graph forms the upper part of the convex hull. Bertoin \cite{B1} studied a non-Gaussian example, establishing the structure of the convex minorant of the Cauchy process. Davydov  \cite{D}  studied these matters for a generic Gaussian process.

We focus here on the set $\mathcal{E}_{X}$ of the extremal points of $C_{X}$, and more precisely on its projection $\Ext_{X}$ onto the horizontal axis, which essentially bears the same properties. Denote by $\mathcal{E}_{X}^+$ and $\mathcal{E}_{X}^-$ the points of $\mathcal{E}_{X}$ that lie respectively in the upper and lower boundary of $C_{X}$, and $\Ext_{X}^+$ and $\Ext_{X}^-$ their respective projections. 
 The points of $\Ext_{\psi}$ are called \emph{extremal times}, and those of $\Ext_{\psi}^+$ (resp. $\Ext_{\psi}^-$), \emph{extremal superior} (resp. \emph{inferior}) \emph{times}. The set $\Ext_{\psi}$ is called the \emph{extremal } set, and $\Ext_{\psi}^+$ (resp. $\Ext_{\psi}^-$) the \emph{extremal superior} (resp. \emph{inferior}) set.
These structures have physical interpretations in a Burgers turbulence with initial data $X$.

Imagine that particles $\{P_{k};\,k\in\mathbb{Z}\}$ are disposed on the real line with respective initial positions $\{k\in\mathbb{Z}\}$ and initial random velocities   $\{v_{k};\,k\in\mathbb{Z}\}$. Assume that every group of particles evolves freely until meeting another group, and when a group of mass $m$ and velocity $v$ collides with a group of mass $m'$ and velocity $v'$, they form a unified clump of particles with mass $m+m'$ and velocity $\frac{mv+m'v'}{m+m'}$, so that mass and momentum are preserved, but not energy. Let $\psi$ be the piece-wise linear function on $\mathbb{R}$, affine on each interval $[k,k+1]$, such that $\psi(k+1)-\psi(k)=-v_{k}$ for all $k$. It is shown in Theorem~\ref{thm:discreteCase} that, if the marginal distributions of the velocities are non-atomic, the sufficient and necessary condition for particle $P_{k}$ to never collide with its left-hand neighbour is that $k$ is in $\Ext_{\psi}^+$.

 In fluid mechanics, this geometrical problem is closely related to the solution of the one-dimensional Burgers equation with vanishing viscosity. Namely, if $v(x,t)$ is the velocity field of an incompressible fluid, Burgers has introduced the equation
\begin{eqnarray}
\label{Burgers}
\partial_{t}v+ \partial_{x}\left(\frac{v^2}{2}\right)=\epsilon \partial^2_{xx} v
\end{eqnarray}
 as a simple model of hydrodynamic turbulence, where $\epsilon$ is the viscosity. Burgers equation is a simplified version of Navier-Stokes equation, and even if it is known among physicists that it does not describe turbulence very accurately, it is broadly used in many physical problems such as shock wave formation in compressible fluids, or formation of large scale structure of the universe. See \cite{W} for more detailed information. 
Frisch explained why a probabilistic description of turbulence is necessary, and there has been since the 80's an abundant literature about the solution of Burgers equation with random  initial conditions  and vanishing viscosity, defined as the limit of the solution as $\epsilon$ tends to $0$.

It is explained in Section~\ref{sec:Physics} how a topological description of the fluid can be deduced from the study of the convex hull of $\psi+f$, where ${\partial_{a}}\psi(.,0)=-v$, and $f$ is a parabolic drift. If the drift is removed, the study of the convex hull of $\psi$ retrieves the asymptotic structure of the fluid when time goes to $\infty$. In particular, the points of $\Ext_{\psi}^+$ are the extremities of shock intervals, where a shock interval contains all the initial positions of particles that end up in the same clump after a finite time. If $a$ is a left and right accumulation point of $\Ext_{\psi}^+$, it means that $a$ has never collided with its  neighbours during the turbulence. Such a point is said to be \emph{Lagrangian regular}. 

Many authors studied the shock structures of fluids fed with random initial data. See \cite{W} for a detailed account of the progress in this direction. The seminal papers \cite{S} and \cite{She}, and \cite{Av95} and \cite{Av95b}, obtained results in a Gaussian framework.  In the field of L\'evy processes, \cite{B4} and \cite{B2} studied the structure of shocks when the initial velocity is a stable L\'evy noise, with emphasis on the set of Lagrangian regular points. 
We are interested mainly in this paper on the asymptotic topology of the shock structure of a turbulence which velocity is a L\'evy noise, and we derive metrical results for a more general class of Markov processes, and for It\^o processes. The assumption of self-similarity present in most precedent works is dropped here, and the main result concerns the general class of L\'evy processes with bounded variation. We will prove in this framework that  with probability $1$ there are no Lagrangian points (i.e. points which never undergo a collision) on a non-compact interval. If the interval is compact, the only particle capable of remaining untouched during the turbulence is the one which initial location is $T$, the point where the initial L\'evy process reaches its maximum. We show that it is indeed a Lagrangian regular point iff $0$ is regular for both half-lines.

For a function $f$ on an interval of $\mathbb{R}$, denote by $\overline{f}$ the smallest concave function greater than $f$, or \emph{concave majorant}, and $\underline{f}=-(\overline{-f})$ its convex minorant. We also have $\mathcal{E}_{\overline{f}}=\mathcal{E}_{f}^+$ and  $\mathcal{E}_{\underline{f}}=\mathcal{E}_{f}^-$. 

We state the results in terms of the convex  hull of a random process $X$.
 It is shown in this paper that, if $X$ belongs to a certain class of Markov processes, called \emph{reversible}, then $\mathbf{E}_{X}$ is a.s. negligible, and so it is for its integrated process. The same result holds if $X$ is an It\^o process. The main result of this paper concerns the case where $X$ is a L\'evy process with bounded variation. The set $\mathbf{E}_{X}^+$ is a.s. countable with accumulation only possible  at the time $T$ where $X$ tends to its maximal value. We also give sufficient and necessary conditions for $T$ to be an accumulation point of $\mathbf{E}_{X}^+$, on its right and on its left, in terms of the regularity of $0$ for the half-line. More generally, if a sufficiently smooth deterministic function $f$ is added up to $X$, we show that the only right (resp. left) accumulation points of $\mathbf{E}_{X}^+$ can occur only at times $a$ where $f'(a)=(\overline{X+f})'(a^+)$ (resp. $f'(a)=(\overline{X+f})'(a^-))$.\\

In the preliminaries, we introduce rigourously all the notation used in the paper, and recall some facts about L\'evy processes which will be useful in the sequel. Then, in Section \ref{sec:Physics}, we explain the relation between the convex hull of a function and the description of the shock structure of a fluid without viscosity. Although it is not an accurate model, the discrete case is presented as well to help comprehend the main aspects of the problem. Finally, we state the results concerning the convex hull for some random processes,  with a complete description for a L\'evy process with bounded variation.\\

 In this paper, we investigate the metrical and topological properties of the set of extremal times $\mathbf{E}_{X}$ of a c\`adl\`ag random process $\{X(a)~;~a\in I\}$, for an interval $I$ of $\mathbb{R}$. The symmetry of the problem allows us  to simply study $\mathbf{E}_{X}^+$, the properties of $\mathbf{E}_{X}^-$ being similar. In this respect, define $X^*(a)=\max(X(a^-),X(a))$, so that $a\in \mathbf{E}_{X}^+$ iff $(a,X^*(a))$ is an extremal point of the graph of $X$. 
 
 For a process $X$ defined on a right neighbourhood of a point $a\in\mathbb{R}$, denote by $R_{a}^+(X)$ the event $$\forall \epsilon>0,\,\exists s\in]a,a+\varepsilon[,\,X(s)\geq X(a)$$ and $R_{a}^-(X)=R_{a}^+(-X)$.  If the process satisfies Blumenthal's \emph{zero-one law}, which is the case for ``nice'' Markov processes, then they both are trivial. Furthermore, for L\'evy processes, the probability of these events do not depend on $a$ if it is a deterministic time, or more generally a stopping time, in this case simply write $R^+(X)$ and $R^-(X)$, and $p^+(X)=\mathbb{P}(R^+(X)),~p^-(X)=\mathbb{P}(R^-(X))$, both of these numbers being in $\{0,1\}$. The value of $p^+(X)$ is related to the regularity of $0$ for the half-line when $X$ has infinite L\'evy measure.
 
 In all the article, for a process $X$ defined on an interval $I$, the notation $M_{a}$ stands for $(a,X^*(a))$, where $a$  belongs to $I$ and might be random.

\section{Preliminaries}
\subsection{L\'evy processes}
\label{sec:LevyPro}

It is known (see \cite{BBook}, Section 1.1) that the sample-paths of a  L\'evy process $Y$ either have bounded variation on every interval a.s., or unbounded variation on every interval a.s. If $Y$ has bounded variation, it  can be decomposed on any interval $I$ as $$Y(a)=b a+X(a),\, a\in I,$$ where the real number $b$ is the \emph{drift}, and $X$  the Poisson component, a pure jump L\'evy process whose characteristic exponent is of the form $$\psi(\theta)=\int_{\mathbb{R}}(\exp(i\theta x)-1)\nu(dx),\theta\in\mathbb{R}$$ where $\nu$ is a measure on $\mathbb{R}$ satisfying $\nu(\{0\})=0$ and $$\int_{\mathbb{R}}(1\wedge |x|)\nu(dx)<\infty.$$
The measure $\nu$ is the \emph{L\'evy measure} of $Y$ and determines the size and frequency of the jumps.
The drift $b$ is irrelevant in our study, because the set $\mathbf{E}_{Y}$ does not depend on $b$, thus we will assume henceforth that $b=0$, i.e $Y=X$ has no drift and only varies by jumps.

Such processes also enjoy the property of being ``flat'' around the origin (see for example \cite{Sk},
Appendix, Th.1).

\begin{Lm}
\label{ComptLocal}
Let $X$ be a L\'evy process with bounded variation and no drift on an interval $I$. Then we have, for any $a\in I$, with probability $1$,
$$\lim_{x\to a}\frac{X(x)-X(a)}{x-a}=0.$$
Furthermore, due to the strong Markov property, the previous relation still holds if $a$ is replaced by a stopping time.
\end{Lm}

There is another relevant feature regarding the extremal set of L\'evy processes. Say that $0$ is \emph{regular for the positive half-line} for $X$ if  $X$ enters the set $]0,\infty[$ at arbitrary small times. Remark that if $X$ has infinite L\'evy measure, it corresponds with our notation to $p^+(X)=1$. Bertoin \cite{B3} established an explicit characterisation of processes $X$ for which $0$ is regular for the positive half-line, in terms of the L\'evy measure $\nu$ of $X$.

\section{Fluid in vanishing viscosity}
\label{sec:Physics}

The study of the convex hull of a function $f$ is connected to the shock structure of a fluid whose initial potential is $f$. It is easier to understand the technical details in a discrete framework, i.e with discrete particles, but the continuous case is more realistic, and is studied in the latter half of this section.

\subsection{Discrete case}
Consider $\{P_{k}~;~k\in \mathbb{Z}\}$, a set of particles disposed on the real line having respective initial positions $\{k~;~k\in \mathbb{Z}\}$. Particle $P_{k}$ has initial random velocity $v_{k}\in \mathbb{R}$ and mass $1$. The rule of evolution is the following: Any clump of particles with velocity $v$ and mass $m$ evolves freely on the line until meeting another clump of velocity $v'$ and mass $m'$, in which case they form a new clump with mass $m+m'$ and velocity $\frac{mv+m'v'}{m+m'}$, so that mass and momentum are preserved.\

  Define the initial potential $\psi$ by $\psi(0)=0$ and $\psi(k+1)-\psi(k)=-v_{k}$. We suppose that the laws of the random variables $v_{k}$ yield that, for all $k<p<r$ in $\mathbb{Z}$, $\frac{\psi(p)-\psi(k)}{p-k}\neq \frac{\psi(r)-\psi(p)}{r-p}$ almost surely. It is the case for instance if any marginal distribution  $(v_{k_1},..,v_{k_{q}})$  on $\mathbb{R}^q,\, q\geq 1$, has no atom. We still call $\psi$ the linear interpolation of $\psi$ on $\mathbb{R}$, obtained by connecting the dots with segments.
\begin{T}
\label{thm:discreteCase}
There is a partition of the set of integers in intervals $\mathbb{Z}=\cup_{q\in \mathbf{E}_{\psi}^+}Z_{q}$, where $Z_{q}$ is  the set of indices of all particles that end up in the same clump as $P_{q}$ after a finite time. Note that  $q\in\Ext_{\psi}^+$ is the left  extremity of $Z_{q}$. In particular, a particle $P_{k}$ never collides with a neighbour  iff $k$ and ${k+1}$ are in $\mathbf{E}_{\psi}^+$.
\end{T}

This theorems enlightens the fact that the study of the concave majorant of $\psi$ brings all the information needed for an asymptotic topological description of the fluid.
\begin{proof}

To understand the behaviour of this particles system when time goes to $\infty$, we investigate conditions on $v_{k}$ so that particle $P_{k}$ hits its left neighbour $P_{k-1}$, for some $k$ in $\mathbb{Z}$. 

This happens if the clump formed by $P_{k}$ and some of its right neighbours, say $\{P_{k},\dots,P_{k+s}\}$, for some $s\in \mathbb{N}^*$, collides with the clump formed by $P_{k-1}$ and some of its left neighbours, say $\{P_{k-1},\dots,P_{k-1-u}\}$ for some $u$ in $\mathbb{N}^*$. This event will occur if the velocity of the right hand clump is smaller than that of the left hand clump, and this writes
\begin{equation*}
\frac{1}{u}(v_{k-1}+\dots+v_{k-u})>\frac{1}{s}(v_{k}+\dots+v_{k+s}).
\end{equation*}
 
Thus, the necessary and sufficient conditions for $P_{k}$ (and its right hand clump) to infinitely avoid a collision with $P_{k-1}$, and maybe other particles $P_{j},j<k$, is

\begin{eqnarray*}
\frac{\psi(k)-\psi(k-u)}{u}>\frac{\psi(k+s+1)-\psi(k)}{s},~s,u\in \mathbb{N}^*,
\end{eqnarray*}
 the equality cases almost never happen, due to the assumptions on the laws of the $v_{k}$.
 This condition exactly means that $k$ is an extremal superior time  of $\psi$. 
 \end{proof}
 
 Considering that $P_{k}$ hits its right neighbour iff $P_{k+1}$ hits its left neighbour, the study of $\mathbf{E}_{\psi}^+$ yields the full description of the shock structure of the system, after sufficiently large times.
 
 \cite{Wi2} described more precisely the shock structure of such a discrete set of particles with random initial velocities driven by Burgers turbulence.

\subsection{Continuous case}
\label{ContCase}

Now, we assume that a fluid with viscosity $\epsilon$ is spread  on an interval $I$ of the real line, and consider that its evolution is ruled by Burgers equation (\ref{Burgers}). Here, $v(a,0)$ denotes the initial velocity of the particle initially located in $a$, and $v(a,t)$ is the  velocity field of the fluid at time $t$ and point $a$. Hopf \cite{H} and Cole \cite{C} derived an explicit solution, stated in terms of the potential
 $\psi$ defined by $\partial_{a}\psi(a,t)=-u(a,t)$ and $\psi(0,0)=0$. The limit of this solution when the viscosity $\epsilon$ tends to $0$ is given, for $x\in \mathbb{R},t\geq 0$, by
\begin{eqnarray}
\label{HopfCole}
\psi(x,t)=\sup_{a \in \mathbb{R}}\left[\psi(a,0)-\frac{(x-a)^2}{2t}\right]=-\frac{x^2}{2t}+\sup_{a \in \mathbb{R}}\left[\psi(a,0)-\frac{a^2}{2t}+\frac{xa}{t}\right].
\end{eqnarray}
Remark the similarity with the Legendre transformation of the function
\begin{eqnarray*}
\psi_t(a)=\psi(a,0)-\frac{a^2}{2t}.
\end{eqnarray*}

Call $a(x,t)$ the greatest value of $a$ for which the supremum in (\ref{HopfCole}) is achieved. The function $x\mapsto a(x,t)$ is  non-decreasing and right-continuous, thus we can define   $x(a,t)$ its right-continuous inverse. Note that the mere existence of this supremum presupposes that $\psi(a,0)=o(a^2)$ when $a$ goes to $\infty$. This is indeed  the case with probability $1$ if $\psi$ is a L\'evy process.

Graphically, if one lets the line with slope $-\frac{x}{t}$ descend vertically until touching the graph of $\psi_{t}$ (assuming that the line was initially high enough), $a(x^-,t)$ is the abscissa of the most leftward contact point between the line and the graph of $\psi_{t}$, while $a(x,t)$ is the most rightward contact point. 
It is clear from this description that for all $x$, $a(x^-,t)$ and $a(x,t)$ are in $\Ext_{\psi_{t}}^+$.

 The application $a \to x(a,t)$ is called \emph{Lagrangian}, and $x\to a(x,t)$ \emph{inverse Lagrangian} function. From the hydrodynamic point of view, $x(a,t)$ is the position at time $t$ of a particle initially located at $a$. If a discontinuity of the inverse Lagrangian occurs at a point $x$: $a(x^-,t)<a(x,t)$, it means that all the particles initially located on $[a(x^-,t),a(x,t)[$ have formed a clump at point $x$ at time $t$, and that is why such an interval is called a \emph{shock interval}. 

Of particular interest are the points $a$ which are in the closure of no shock interval at time $t$. Such points are called \emph{Lagrangian regular points}, and are characterised as the points of $\Ext_{\psi_t}^+$ that are left and right accumulation points of $\Ext_{\psi_t}^+$. Denote their set by $\Lr_{\psi_t}$. The Lagrangian regular points represent the initial locations of particles that have not been involved in any shock up to time $t$.

 Thus, if we call $\Sh_{\psi_t}=\{a(x^-,t),a(x,t):\,a(x^-,t)<a(x,t)\}$ the extremities of shock intervals,  we have the disjoint union $\Ext_{\psi_t}^+=\Sh_{\psi_t}\cup \Lr_{\psi_t}$. 
 Hence the topology of $\mathbf{E}_{\psi_t}^+$ contains the whole description of the shock structure of the fluid at time $t$. Notice the similarity of this result with theorem \ref{thm:discreteCase}.\\

In conclusion, in order to obtain a complete description of the fluid  whose initial velocity is a L\'evy noise, our aim is to study $\mathbf{E}_{\psi+f}^+$, if $\psi$ is the realisation of a L\'evy process and $f$ is the smooth drift $f(a)=-\frac{a^2}{2t}$.  We conduct the study for general smooth $f$, and we will consider the special case $f(a)=-\frac{a^2}{2t}$ in Section~\ref{sec:burgers-results}.

Since the drift goes to $0$ as time goes to $\infty$, it is natural to wonder about the asymptotic behaviour of $\Ext_{\psi_{t}}$. The following result yields a form of continuity for $t\to \mathbf{E}^+_{\psi_{t}}$.

Denote by $B(0,\varepsilon)$ the open ball with center $0$ and radius $\varepsilon$ in $\mathbb{R}$, so that for a set $M\subseteq \mathbb{R}$, $M+B(0,\varepsilon)$ is the $\varepsilon$-neighbourhood of $M$.

\begin{Ps}
\label{UppSemCon}
Let $g$ be a function on a bounded interval $I$ of $\mathbb{R}$, and $f$ a continuous function on $I$. For a real number $u$, call $g_{u}=g+uf$ and $\mathbf{E}^+_{u}=\mathbf{E}^+_{g_{u}}$. Assume that 
\begin{eqnarray}
\label{eq:clear}
{\rm{cl}}({\rm{gr}}(g_{u}))\cap {\rm gr}(\overline{g_{u}})=\mathcal{E}_{g_{u}},
\end{eqnarray}
i.e the graph of $g_{u}$ only approaches the boundary of its convex hull in its extremal points.
Then,
\begin{eqnarray*}
d_{\mathcal{H}}(\mathbf{E}^+_{u+h},\mathbf{E}^+_{u}) \xrightarrow[h \to 0]{}0,
\end{eqnarray*}
where $d_{\mathcal{H}}$ denotes the Hausdorff distance between closed sets, defined by
\begin{eqnarray*}
d_{\mathcal{H}}(L,M)=\inf\{r>0~;~L \subseteq( M+B(0,r)), M \subseteq (L+B(0,r))\}.
\end{eqnarray*}

\end{Ps}

\begin{proof}
Let $\epsilon>0$. We need to find $\eta>0$ such that $\mathbf{E}^+_{u-h} \subseteq \mathbf{E}^+_{u}+B(0,\epsilon)$ for $0<h<\eta$. Since we can replace $g$ by $-g$, it will be sufficient.

 Since $\mathbf{E}^+_{u}$ is closed,  $I\setminus \mathbf{E}^+_{u}$ can be decomposed into a countable set of disjoint open intervals $(I_{n,u})_{n \in \mathbb{N}}$.
Put $K_{u,\epsilon}=I\setminus (\mathbf{E}^+_{u}+B(0,\epsilon))$. It intersects a finite number of $I_{n,u}$, the ones whose size is greater than $\epsilon$. Call ${J}_{1},...,{J}_{q}$ the closed intervals such that $K_{u,\epsilon}=\cup_{i=1}^q {J}_{i}$. 

Take $i$ in $\{1,...,q\}$.  Since by hypothesis the closure of the graph of $g_{u}$ does not touch that of $\overline{g_{u}}$ outside an extremal point, the graphs of $g_{u}$ and $\overline{g_{u}}$ have positive distance on each of the ${J}_{i},1\leq i \leq q$. In consequence, there is $h_{i}>0$ such that the graphs of $g_{u-h_{i}}$ and $\overline{g_{u-h_{i}}}$ are disjoints  above ${J}_{i}$. In particular, $\mathbf{E}^+_{u-h_{i}} \cap {J}_{i}=\emptyset$. Taking $h=\min_{i}(h_{i})>0$ yields the result.
 \end{proof}

For typical realisations of the irregular random processes studied in this paper, (\ref{eq:clear}) will be satisfied almost surely for almost all positive $u$. But also, with probability one, there will be some (random) $u$ such that it is not satisfied. 

The asymptotic case corresponds to the value $t=+\infty$, or $u=0$. The previous theorem allows us to state that the structure of the fluid at sufficiently large times tends to the extremal set of $\psi(.,0)$, without drift. This limiting set $\mathbf{E}_{\psi}^+$ is studied in detail in the next section.

\section{Extremal set of random processes}

The first part of the results consists in proving that the negligibility of the extremal set of L\'evy processes and It\^o processes occurs almost surely, if the underlying measure is Lebesgue measure $\lambda$. Then the topology of the extremal set of L\'evy processes with bounded variation is more profoundly investigated, and stronger results are obtained. 

\subsection{Negligibility results}

To establish the negligibility, the main tool is Fubini's theorem, which enables us to use the following lemma.
In all this section, $I$ stands for an interval of $\mathbb{R}$.
\begin{Lm}
\label{LmNgl}
Let $X$ be a separable process on $I$. If
\begin{eqnarray*}
\A a \in I, ~~\mathbb{P}(a \in \mathbf{E}_{X}^+)=0,
\end{eqnarray*}
then we have a.s. $\lambda(\mathbf{E}_{X}^+)=0$.
\end{Lm}
\begin{proof}
The application $(\omega,a) \to {1}_{(a ~\in~ \mathbf{E}_{X}^+)}$ is measurable:
\begin{multline*}
\left\{(\omega,a)~;~a \notin \mathbf{E}_{X}^+(\omega)\right\}  =\\ 
\left\{(\omega,a)~;~X(\omega,a) \leq  \limsup_{{s,v ~\in \mathbb{Q}^2}\atop {s<a<v}} (a-s)X(\omega, v)+(v-a)X(\omega, s)\right\} .
\end{multline*}
Fubini's theorem gives us
\begin{eqnarray}
\label{Fubini}
0=\int_{I}\mathbb{P}(a \in \mathbf{E}_{X}^+)~ dt = \int_{\Omega} \lambda(\{a~;~a\in \mathbf{E}_{X}^+\})~  \mathbb{P}(d\omega)
\end{eqnarray}
and the proof is complete.
\end{proof}

Due to this lemma and the nice homogeneity properties of L\'evy processes, we can establish the negligibility of the extremal set for any L\'evy process, but the space homogeneity is actually not fully required, and  the result holds under weaker hypotheses. Besides the strong Markov property, we require from a process $X$  the triviality of the algebra $\bigcap_{s \downarrow 0}\sigma(X_{s})$, which is called Blumenthal zero-one law. It is also required that, a.s., the process is a.e. continuous.
\begin{D}
Let $X$ be a Markov process. $X$ is said to be $\emph{reversible}$ if, for any $a\in I$, the processes $\widehat{X}(s)=X(a+s)-X(a)$ and $\widecheck{X}(s)=X(a)-X((a-s)^-)$ (for $s$ such that both expressions make sense) are Markov, and have the same distribution.
\end{D}
For example, L\'evy processes are reversible.
\begin{T}
\label{ThmNegRev}
Let $X$ be a reversible Markov process. Then 

\begin{description}
\item[ (i)] $\lambda(\mathbf{E}_{X})=0$ a.s.
\item[ (ii)] Let $Z$ be a process which derivative is $X$. We have a.s. $\lambda(\mathbf{E}_{Z})=0$.
\end{description}

\end{T}

This result is proved in paragraph \ref{PrfNegRev}. Briefly, the idea is to use Lemma \ref{LmNgl} and so we have to prove that any given $a$ is almost never in the extremal set. The property of being in the extremal set depends geometrically  on the past and the future of $a$, and the reversibility property ensures us that the  behaviour of past and future is symmetric. The zero-one law for the process $\widehat{X}-\widecheck{X}$ finishes the proof.\\

The following result deals with It\^o processes. We consider hereafter $I=\mathbb{R}_+$ for the sake of simpler statements.
 Let $({B}_{a};~a \geq 0)$ be a standard brownian motion and $\mathcal{F}=(\mathcal{F}_{a};~a \geq 0)$ its natural filtration.

\begin{D}
Let $\mathcal{M}_{2}^{loc}$ be the class of  $\mathcal{F}$-adapted processes $\psi$ on $\mathbb{R}^+$ that satisfy

\begin{eqnarray*}
\A a>0,~  \mathbb{E}\left(\int_{0}^a \psi^2(s) d s \right)<\infty.
\end{eqnarray*}
\end{D}

\begin{T}
\label{ItoNeg}
Let $X$ be an It\^o process. $X$ can be written under the general form
\begin{eqnarray}
\label{EquIto}
X(a)=\int_{0}^a \psi(s) d B(s) + \int_{0}^a \phi(s) d s,
\end{eqnarray}
where  $\psi,\phi \in \mathcal{M}_{2}^{loc}$.

Denote by $Y(a)$ the first term in (\ref{EquIto}) and $Z(a)$ the  second one. Assume that, with probability one, there is no non-empty open interval where  $Y$ vanishes. Then, a.s. $ \lambda(\mathbf{E}_{X})=0$.
\end{T}

The hypothesis on $Y$ is necessary because if $Y=0$, $Z$ can be chosen to be any smooth function, bringing a heavy extremal set. In other words, the irregularity of $Y$ is sufficiently strong to overbalance the regularity of $Z$. In the proof (paragraph \ref{PrfNegIto}), we use Lemma \ref{LmNgl} again, showing that for any $a$, $Y$ is not locally Lipschitz in $a$ while $Z$ is, so that the large fluctuations of $Y$ around $a$ make it a.s. impossible for $a$ to be in the extremal set.

\subsection{L\'evy processes with bounded variation}

\label{plvbsd}
 L\'evy processes are often considered as initial data for Burgers turbulence. We focus here on the class of L\'evy processes with bounded variation, and also assume that the drift is null. In this case  Lemma \ref{ComptLocal} indicates  the local behaviour of the process after zero, and consequently after each stopping time. Since the process is pure-jump, and its jumps are all simultaneously stopping times, it is easier to apprehend the structure of the graph with Markovian techniques.

In this section, $X$ is a L\'evy process with bounded variation and no drift, and $f$ is a smooth function. Hypotheses on $f$ will be made more precise later. Let us start with the case where $f$ is convex.\\

\begin{Lm}
\label{Convex}
Let $g$ and $f$ be two functions on $I$, $f$ being convex. Let $a$ be an extremal superior time of $g+f$, then $a$ is an extremal superior time of $g$.
\end{Lm}

\begin{proof}
Let us pick $s$ that is not extremal superior for $g$. We will show that $s$ is not extremal superior for $g+f$ either. We can find $u<s<v$ and $\alpha, \beta=1-\alpha$ in $]0,1[$ such that
\begin{eqnarray*}
s=\alpha u+\beta v,\\
g(s)\leq\alpha g(u)+\beta g(v).
\end{eqnarray*}
The convexity of  $f$ gives us
\begin{eqnarray*}
f(s)\leq \alpha f(u) + \beta f(v),
\end{eqnarray*}
 and by adding up we have
\begin{eqnarray*}
(g+f)(s) \leq \alpha (g+f)(u)+ \beta(g+f)(v),
\end{eqnarray*}
which proves that $s$ is not extremal superior for $g+f$.
\end{proof}

\begin{C}
\label{cor:cvx}
If $f$ is convex, $\mathbf{E}_{X+f}^+ \subseteq \mathbf{E}_{X}^+ \subseteq \mathbf{E}_{X-f}^+.$
\end{C}
Adding up a convex function to $X$ thins its extremal superior set, although adding a concave function widens it.

 The next theorem sums up the main result obtained about concave functions. 
\begin{T}
\label{EnvProLevDri}
Suppose that $f$ is a concave function of class $\mathcal{C}^1$.
Put $Y=X+f$.

 Then, for all $a \in \mathbf{E}_{Y}^+,~ a$ is left isolated (resp. right isolated) in $\mathbf{E}_{Y}^+$ if $\overline{Y}'(a^-)\neq f'(a)$ (resp. $\overline{Y}'(a)\neq f'(a)$).
\end{T}
The proof is in Section \ref{PrfEnvProLevDri}. The  accumulation points of $\mathbf{E}_{Y}^+$, if there are any, are the points where the derivatives of $\overline{Y}$ and $f$ exactly coincide.
With the help of Corollary \ref{cor:cvx}, this theorem can easily be generalised to functions $f$ for which the interval $I$ can be decomposed in a countable family of intervals above which $f$ is either concave or convex.

Applying Theorem \ref{EnvProLevDri} with $f=0$ yields the asymptotic evolution of a fluid, as explained in Proposition~\ref{UppSemCon}
. In this case, the only point where it is possible to have $\overline{Y}'(a)=f'(a)=0$ is the supremum of $X$. The topology of the convex hull around the maximal value then depends on the \emph{regularity of the half-line} for $X$, see Section \ref{sec:LevyPro}.
The following theorem gives an exhaustive topological description of the extremal set.

Remark that if the L\'evy measure is finite, $X$ is a.s. piece-wise constant, and $\mathbf{E}_{X}$ is discrete in $\mathbb{R}^+$. For the sake of more simple statements, we assume in the sequel that the L\'evy measure is infinite.

\begin{T}
\label{EnvLevyVB}
Let $X$ be a bounded variation L\'evy process with infinite L\'evy measure and no drift on a compact interval $I$ of $\mathbb{R}$. Let $T$ be the time, a.s. unique,  that satisfies $X^*(T)=\sup_{a\in I}X(a)$. Then a.s. $\mathbf{E}_{X}^+$ contains only $T$ and some jump times of $X$, and its unique accumulation point is  $T$.

The time $T$ is right isolated in $\mathbf{E}_{X}^+$ iff $0$ is irregular for the positive half-line, and left isolated iff $0$ is irregular for the negative half-line.
\end{T}

\indent This theorem gives us a fairly good understanding of the shape of $X's$ concave majorant and extremal set. We do not have a complete quantitative description as \cite{G} and \cite{P} got for brownian motion, but we can fully understand the topological structure of the set. 

The example of the subordinators perfectly fit this framework. If $X$ is a subordinator on $[0,1]$ with infinite L\'evy measure, then $1$ is the only accumulation point of $\Ext_{X}^+$.

\subsection{Shock structure of a Burgers turbulence}
\label{sec:burgers-results}

We consider here a fluid on an interval $I$,   ruled by Burgers equation (\ref{Burgers}), with initial potential $\psi$. The considerations of Section \ref{ContCase} along with the results obtained in the previous section yield the following results.

\begin{T}
Assume that $\psi$ is a L\'evy process with bounded variation, no drift and infinite L\'evy measure, and $I$ an interval of $\mathbb{R}$.  For all random or deterministic time $a$, call $P_{a}$ the particle whose initial position is $a$.

Assume that $I$ is compact. Let $T$ be the almost unique point of $I$ where $\psi^*$ reaches its maximum. 
\begin{description}
\item[(i)] The only particle that might  never be involved in any shock is the particle $P_{T}$.
\item[(ii)] The particle $P_{T}$ does not undergo any leftward (resp. rightward) collision iff $0$ is regular for the positive half-line (resp. negative half-line).
\item[(iii)] The set $\Ext_{\psi}^+$ is countable and discrete away from $T$ and only contains jump times, apart from $T$.
\item[(iv)] For all $\varepsilon>0$, there is a time $t_{\varepsilon}$ after which each particle $P_{a}$ that is either the extremity of a shock interval or a Lagrangian regular point is at distance at most $\varepsilon$ from $\Ext_{\psi}^+$.
\item[(v)] At an arbitrary time $t>0$, any Lagrangian regular point $a$ at time $t$ must satisfy $\overline{\psi}'_{t}(a)= \overline{\psi}'_{t}(a^-)= -\frac{a}{t}$.
\end{description}

We have the following consequences  on the half-line $I=\mathbb{R}_{+}$. Note that \emph{\textbf{(v)}} still holds in this case.

\begin{description}
\item[(i)'] Every particle is involved in a shock in a finite time.
\item[(ii)'] The set $\Ext_{\psi}^+$ is discrete and only consists of positive jump times.
\item[(iii)'] The set of shock points and Lagrangian regular points tends to locally shrink around $\Ext_{\psi}^+$, in the sense of \emph{\textbf{(iv)}}.

\end{description}

\end{T}

\begin{proof}
{\bf (i),(ii),(iii)} are direct consequences of Theorem~\ref{EnvLevyVB}, and \textbf{(iv)} comes from Proposition~\ref{UppSemCon}.

{\bf (v)}  is a consequence of Theorem~\ref{EnvProLevDri}. 

{\bf (i)',(ii)'}: According to \textbf{(i)}, an accumulation point of $\Ext_{\psi}^+$  is supposed to achieve the maximum of $\psi^*$ on any compact interval, thus it also supposed to achieve the supremum of $\psi^*$ on $I$. Since a non-null L\'evy process without drift  is not majorised, there is no such accumulation point.

{\bf (iii)'} comes from \textbf{(iv)}.
\end{proof}

These results complete those of \cite{B2}, who studied the shock structure at a given time $t>0$ when the initial L\'evy process $\psi$ is a stable L\'evy noise with index $\alpha>\frac{1}{2}$. Note that such processes only have bounded variation if $\alpha<1$ as well. \cite{Win} also established asymptotic results, showing that the shock structure can be thought of as a point process in $\mathbb{R}_{+}\times \mathbb{R}$, where a point $(m,v)$ represents a clump of mass $m$ and velocity $v$, when the initial velocity is a self-similar Markov process.

By combining \textbf{(ii)'} and \textbf{(iii)'}, one can attach to each point $a$ in the discrete set $\Ext_{\psi}^+$ a function $\{\varepsilon_{a}(t);\,t>0\}$, that vanishes when $t\to\infty$, such that, at any time $t>0$, the set $\Ext_{\psi_{t}}^+$of shock points and Lagrangian regular points is contained in $\cup_{a\in\Ext_{\psi}^+}[a-\varepsilon_{a}(t),a+\varepsilon_{a}(t)]$. 

For instance, this holds if $\psi$ is a subordinator. In this case the shock structure is asymptotically discrete. This agrees with the results of Bertoin \cite{B2}, who stated that the shock structure is already discrete at any finite time $t>0$ when $\psi$ is a stable subordinator of index $\alpha\in ]\frac{1}{2},1[$. Thus in this particular situation, the shock structure does not drastically change as time goes to $\infty$. It would be interesting to see if the latter statement holds without the hypothesis of self-similarity.

\section{Proofs}
We introduce the temporal translation operator $\theta$:

\begin{D}
Let $T\in \mathbb{R}^+$, and $Y$ be  a function on an interval $I$ taking values in $\mathbb{R}$ or $\mathbb{R}^2$. Then we set
\begin{eqnarray*}
\theta_{T}(Y)(a)=Y(T+a)
\end{eqnarray*}
for $a$ and $T$ such that  $a,T+a\in I$.
\end{D}

\subsection{Reversible Markov processes}
\label{PrfNegRev}

\begin{proof}[Proof of Theorem \ref{ThmNegRev}]
 Let $a~\in I$. We have to show that $a$ is almost never an extremal superior time.

We set $M_{a}=(a,X^*(a))$. By hypothesis, for $a$ in a set of full measure,  with probability one,  $X^*({a})=X(a)$ (because $X$ is continuous in $a$), and we only consider sample-paths  for which $X^*({a})=a$ in the sequel.

We know that $a$ is not an extremal superior time as soon as we can find $s>0$ such that
\begin{eqnarray}
\label{Cond}
X(a-s)+X(a+s)\geq 2X^*(a),
\end{eqnarray}
because then $M_{a}$ would be under the segment $[M_{a-s},M_{a+s}]$. 
We will show that we can a.s. find $s$ such that (\ref{Cond}) is satisfied.

We set, for $s\geq 0$,
\begin{eqnarray*}
\widehat{X}(s)=X(a+s)-X(a),~\widecheck{X}(s)=X(a)-X((a-s)^-)
\end{eqnarray*} and
\begin{eqnarray*}
Y(s)=\widehat{X}(s)-\widecheck{X}(s)=X((a-s)^-)+X(a+s)-2X(a).
\end{eqnarray*}
By reversibility of $X$, $\widecheck{X}$ and $\widehat{X}$ have the same law and are Markov processes for $s$ small enough. 

By the strong Markov property, for $x \in \mathbb{R}$, the conditional processes $\widehat{X}_{x}=(\widehat{X}|X(a)=x)$ and $\widecheck{X}_{x}=(\widecheck{X}|X(a)=x)$ are independent.
As a consequence, $Y_{x}=\widehat{X}_{x}-\widecheck{X}_{x}$ is Markov and symmetric ($-Y_{x} \stackrel{(d)}{=}  Y_{x}$). Then we have $\mathbb{P}(R_{0}^+(Y_{x}))=\mathbb{P}(R_{0}^-(Y_{x}))\in\{0,1\}$. Both probabilities cannot be equal to $0$, thus they are equal to $1$.

We have  $p^+(Y)=\int_{\mathbb{R}}\mathbb{P}(R_{0}^+(Y_{x}))\mathbb{P}(X(a)\in d x)=1$ and so we can a.s. find $s$ arbitrarily close to $0$ such that $Y(s)\geq 0$, which proves that $X$ satisfies (\ref{Cond}).
So, a.s. $a \notin \mathbf{E}_{X}^+$.
Thanks to Lemma \ref{LmNgl} we can conclude that a.s. $\lambda(\mathbf{E}_{X}^+)=0$.

Given that $-X$ satisfies the same hypotheses, a.s. $\lambda(\mathbf{E}_{-X}^+)=\lambda(\mathbf{E}_{X}^-)=0$ and a.s. $\lambda(\mathbf{E}_{X})=0$, which proves \textbf{(i)}.\\

For showing \textbf{(ii)}, we assume without loss of generality that $[0,a]$ is in the interior of $I$. Let us now consider $Z(a)=\int_{0}^a X(u) d u$. We have to show, like previously, that, with probability $1$, there is $s>0$ such that
\begin{eqnarray*}
\int_{a-s}^a X(u) d u \geq \int_{a}^{a+s} X(u) d u .
\end{eqnarray*}
We set, for $s$ small enough,
\begin{eqnarray*}
\widehat{Z}(s)=\int_{a}^{a+s} X(u) d u= \int_{0}^s \widehat{X}(u) d u\text{~ and ~}\widecheck{Z}(s)=\int_{a-s}^{a} X(u) d u = \int_{0}^s \widecheck{X}(u) d u.
\end{eqnarray*}
Put also  $W(s)=\widehat{Z}(s)-\widecheck{Z}(s)=\int_{0}^sY(s)d s$.

For $x \in \mathbb{R}$, we index by ``$x$'' the variables conditioned by $X(a)=x$.
The process $Z$ is not Markov, but $(Z,X)$ is Markov, and so $(W_{x},Y_{x})$ is Markov too, and has a symmetric law, as the subtraction of two independent Markov processes with the same law.
The \emph{zero-one law} ensures us again that a.s. $W$ is negative arbitrarily close  to $0$, and (\ref{Cond}) is satisfied for $Z$. So a.s. $a \notin \mathbf{E}_{Z}^+$.
 Thanks to Lemma \ref{LmNgl} we can conclude that a.s. $\lambda(\mathbf{E}_{Z}^+)=0$. We arrive at $\lambda(\mathbf{E}_{Z})=0$ a.s.
\end{proof}

\subsection{It\^o processes}
\label{PrfNegIto}

\begin{proof}[Proof of Theorem \ref{ItoNeg}]

$X$ is written under the form
\begin{eqnarray}
\label{Ito}
X(a)=\int_{0}^a \phi(s) d B(s) + \int_{0}^a \psi(s) d s
\end{eqnarray}
 where $(B,\mathcal{F})$ is a standard brownian motion with its filtration and $\phi,\psi \in \mathcal{M}_{2}^{loc}$.\\
Denote by $Y(a)$ the first term in (\ref{Ito}) and $Z(a)$ the second one.
Recall that, by hypothesis, almost surely, $Y$ is null on no interval. 

A preliminary remark is that, given $\tau$ continuous non-decreasing on $\mathbb{R}^+$, $R_{\tau(a)}^+(X)$ is realised if  $X\circ \tau(a+v)-X\circ \tau(a)>0$  for arbitrarily small values of $v$. Moreover, $\tau$ can be random. An ideal candidate for this random time change is the function given by the following theorem: (See \cite{Y}, ch.V.1)
\begin{T}(Dubins-Schwarz)\\
Let $M$ be a local continuous martingale for $\mathcal{F}$. Let $\{\langle M \rangle_{a};\, a\in I\}$ be its quadratic variation. Define 
\begin{eqnarray*}
\tau(a)=\inf\{s ;\,  \langle M \rangle_{s}\geq a\}.
\end{eqnarray*}
  Then $
W=M\circ \tau $ is a standard brownian motion for the filtration $\{\mathcal{F}_{\tau(a)};\,a\in I\}$.
\end{T}

This theorem applies to $Y$ which is a local continuous martingale. So that the time change of $Y$ is continuous, we need that $a \to\langle Y \rangle_{a} $ is constant on no interval. The process $Y$ being a local continuous martingale, on any deterministic interval where $\langle Y \rangle_{a} $ is constant, $Y$ has bounded variation and is hence null. Thus the hypotheses imply that $\tau$ is continuous.

We apply the time change  to $Y$ and keep the same notation.
$$X\circ \tau(a)=W(a)+Z\circ \tau(a).$$

We know that, for any positive $a$, $\limsup_{v \downarrow 0} \frac{W(a+v)-W(a)}{v}=\infty$ a.s.. So, if $Z\circ \tau$ is locally Lipschitz, its contribution in the increasing rate is negligible compared to that of $W$, and
\begin{eqnarray*}
\limsup _{v \downarrow 0} \frac{X\circ \tau(a+v)-X\circ \tau(a)}{v}=\infty.
\end{eqnarray*}
 We finally have 
 $$\mathbb{P}((R_{\tau(a)}^+(X))^c) \leq \mathbb{P}(a \in L_{\tau}^c),$$ where $L_{\tau}$ is the set of points where $\tau$ is Lipschitz. ($Z$ is a.s. an absolutely continuous function, hence  Lipschitz everywhere).
Using Fubini yields
\begin{eqnarray*}
\int_{\Omega} \lambda(\{a~;(R_{\tau(a)}^+(X))^c\})  \mathbb{P}(d\omega) \leq \int_{0}^\infty \ \mathbb{P}(a \in L_{\tau}^c) d a \leq \int_{\Omega} \lambda(L_{\tau}^c)  \mathbb{P}(d\omega).
\end{eqnarray*}
We need a  well known result.
\begin{Lm}[Riesz-Nagy]
If $f$ is a non-decreasing real function, it is differentiable almost everywhere. 
\end{Lm}
Hence, $\tau$ is locally Lipschitz a.e. and a.s. $\lambda(L_{\tau}^c)=0$, which implies that for almost all $ a$,  $R_{\tau(a)}^+(X)$ occurs a.s.  Since $\tau$ is continuous and unbounded, $\{\tau(a):\,a\in\mathbb{R}_{+}\}=\mathbb{R}_{+}$, thus for almost all $a\in \mathbb{R}_{+}$, $R_{a}^+(X)$ occurs a.s. It means that $X$ takes values larger than $X(a)$ arbitrarily close from $a$.

As $X, -X, \widecheck{X}$ and $-\widecheck{X}$ satisfy the same hypotheses, they also  realise $R_{a}^+$ with probability $1$ in every $a$. We  then say that for any $a$, a.s. $M_{a}=(a,X^*(a))$  is strictly included in the convex polygon $M_{a-s_{1}}M_{a-u_{1}}M_{a+s_{2}}M_{a+u_{2}}$ for some $u_{1},u_{2},s_{1},s_{2}>0$, and so is not in $\mathcal{E}_{X}$.

\end{proof}

\subsection{L\'evy processes with smooth drift}
\label{PrfEnvProLevDri}
For the proof of Theorem \ref{EnvProLevDri}, $X$ denotes a L\'evy process with bounded variation and no drift  on $I$, $f$ is a concave function of class $\mathcal{C}^1$  on $I$, and $Y=X+f$. 

Without loss of generality, we only consider the case where $I=[0,+\infty[$, and $X(0)=f(0)=0$.

The proof is based on the use of some specific stopping times,  defined in the following.
\begin{D}
For $\mu>0$, we set $$S_{1,\mu}(X,f)=\inf\{a>0~;~Y(a)-Y(0)> (f'(0)+\mu) a\},$$ and $S_{k+1,\mu}=S_{1,\mu}\circ \theta_{S_{k,\mu}}+S_{k,\mu},\, k\in \mathbb{N}^*$.

For $u \in \mathbb{Q}\cap I$, we set $S_{k,\mu,u}=S_{k,\mu}\circ \theta_{u}$.\end{D}

We write for short $S_{k,\mu}$ in the following, but one has to keep in mind that the temporal translation is applied to both function $X$ and $f$, so that in the definition of $S_{k+1,\mu}$ one compares the increment of $X$ with the derivative of $f$ in $S_{k,\mu}$, and not in $0$.

Those times  are called \emph{exceeding times}, because they correspond to moments where $Y(a+s)$ exceeds the line with slope $f'(a)+\mu$ passing through $(a,Y(a))$, for some $a$.
\begin{Ps}
For each $\mu>0, u\in \mathbb{Q}\cap I $, $\{S_{k,\mu,u};\,k\geq 1\}$ is a sequence of stopping times that goes to $\infty$.
\end{Ps}
\begin{proof}
We are going to show it only for $u=0$, the general case can be easily deduced.

The stopping time aspect doesn't raise any problem.
We define
\begin{eqnarray*}H_{k+1,\mu}=\inf\{a>S_{k,\mu}~;~X(a)-X(S_{k,\mu})> \mu( a-S_{k,\mu})\}.
\end{eqnarray*}
By concavity, 
\begin{eqnarray}
\label{eq1}
& f(S_{k+1,\mu})-f(S_{k,\mu}) \leq f'(S_{k,\mu})(S_{k+1,\mu}-S_{k,\mu}).
\end{eqnarray}
We also have, by the definition of $S_{k+1,\mu}$,
\begin{multline}
\label{eq2}
  X(S_{k+1,\mu})+f(S_{k+1,\mu})-X(S_{k,\mu})-  f(S_{k,\mu})\\
\geq  (f'(S_{k,\mu})+\mu)(S_{k+1,\mu}-S_{k,\mu}).
\end{multline}
So, the subtraction (\ref{eq2})-(\ref{eq1}) yields,
\begin{multline*}
X(S_{k+1,\mu})-X(S_{k,\mu}) \geq  (f'(S_{k,\mu})+\mu)(S_{k+1,\mu}-S_{k,\mu}) -f'(S_{k,\mu})(S_{k+1,\mu}-S_{k,\mu})\\
 \geq\mu (S_{k+1,\mu}-S_{k,\mu}),
\end{multline*}
and so $S_{k+1,\mu} \geq H_{k+1,\mu}$ by the definition of $H_{k+1,\mu}$. Moreover, $H_{k+1,\mu}-S_{k,\mu}$ has the same law as  $\inf\{a>0;\,X(a)-X(0)\geq \mu a\}$. By independence of the increments,  $\{H_{k+1,\mu}-S_{k,\mu};\,k\geq 1\}$ is a sequence of iid random variables strictly positive, their sum tends a.s. to infinity, and so does the sum of the $\{S_{k+1,\mu}-S_{k,\mu};\,k\geq 1\}$ since $S_{k+1,\mu}-S_{k,\mu}\geq H_{k+1,\mu}-S_{k,\mu}$.
\end{proof}

We call $\mathbf{J}_{X}^+$ the set of positive jump times of $X$, containing $0$ and $\infty$ by convention.\\

\begin{Ps}
\label{PpsPextJump}
We can find a countable set $\Lambda$ dense in $\mathbb{R}^*_{+}$ such that with probability $1$  $$\{S_{k,\mu,u};\,u \in \mathbb{Q}\cap I, k \in \mathbb{N}, \mu \in \Lambda\} \subseteq \mathbf{J}_{X}^+.$$
\end{Ps}

\begin{proof}
Without loss of generality, we only treat the case $k=1,u=0$. If we show that for all $\mu>0$
\begin{equation}
\label{eq:polar}
\mathbb{P}(S_{1,\mu}\in \textbf{J}_{X}^+)=1,
\end{equation}
Fubini's theorem will give us the result. With the vocabulary of \cite{BBook}, (\ref{eq:polar}) is equivalent to the polarity of the set $\{0\}$ for the process $\{X(a)-\mu a;\,a\geq 0\}$. According to \cite{BBook}, Pr.~2-(ii), the real part $\mathfrak{R}(\varphi(x))$ of $X$'s L\'evy characteristic exponent $\varphi$ is in $o(|x|)$ when $|x|\to\infty$, thus \cite{BBook},Th.~16-(ii) implies that $\{0\}$ is essentially polar, hence polar (see the comments of the theorem).
\end{proof}

Up to restricting the universe  to a smaller set of full measure, we suppose from now on that for all $\mu\in\Lambda,u\in\mathbb{Q}\cap I$, $S_{1,\mu,u}\in \mathbf{J}_X^+$.

\begin{Ps}
\label{PpsJumpPext}
Let $a$ be an extremal superior time such that $f'(a)<\overline{Y}'(a^-)$.
Then we can find $u \in \mathbb{Q}\cap I, k \in \mathbb{N}, \mu \in \Lambda$ such that $a=S_{k,\mu,u}$. In particular $a$ is a positive jump time.
\end{Ps}

\begin{proof}
We choose $\mu \in \Lambda$ such that $f'(a)<f'(a)+\mu<\overline{Y}'(a^-)$ and $u \in \mathbb{Q}\cap ]0,a[ $ such that for each $s \in [u,a]$ we have $f'(s)+\mu<\overline{Y}'(a^-)$. It is possible because $f'$ is continuous.

Pick $k$ such that $S_{k,\mu,u}\leq a< S_{k+1,\mu,u}$. 
If $S_{k,\mu,u}<a$, then, since $a<S_{k+1,\mu,u}$,
\begin{eqnarray*}
{Y}^*(a)<{Y}(S_{k,\mu,u})+(a-S_{k,\mu,u})(f'(S_{k,\mu,u})+\mu).
\end{eqnarray*}
The concavity of $\overline{Y}$ yields
\begin{eqnarray*}
\overline{Y}(a) \geq \overline{Y}(S_{k,\mu,u}) + (a- S_{k,\mu,u})\overline{Y}'(a^-) \geq Y(S_{k,\mu,u})+(a-S_{k,\mu,u})(f'(S_{k,\mu,u})+\mu).
\end{eqnarray*}
The point $a$ being extremal superior, $\overline{Y}(a)=Y^*(a)$, thus the two last inequalities are contradictory, and we have $S_{k,\mu,u}=a$.
\end{proof}

\begin{proof}[Proof of Theorem \ref{EnvProLevDri}]
Let $a \in \mathbf{E}_{X}^+$ such that $f'(a) \neq \overline{Y}'(a^-)$.\\
If $f'(a) < \overline{Y}'(a^-)$, we have just seen that $a \in \mathbf{J}_{X}^+$, thus it is left isolated in $\mathbf{E}_{X}^+$.

If $f'(a) > \overline{Y}'(a^-)$, applying the same logic to the reversed process $\widecheck{Y}$ and its convex minorant, it is clear that $a$ is a positive  jump time for $\widecheck{Y}$, thus a negative jump time for $Y$. Since the set of all jump times is a countable set of stopping times, by Lemma  \ref{ComptLocal} we have $$\lim_{s \downarrow 0} \frac{Y(b-s)-Y(b^-)}{-s} = f'(b)$$ simultaneously for all jump times $b$. By hypothesis $f'(a) > \overline{Y}'(a^-)$, so the graph of $Y$ on the left of $a$ is located under the line passing through $M_{a}$ with slope $\overline{Y}'(a^-)$, and $a$ is left isolated in $\mathbf{E}_{Y}^+$.

We can apply a similar argument on the right of $a$.\\
\end{proof}

\subsection{L\'evy processes with no drift}

Here, $X$ is a L\'evy process with infinite L\'evy measure, bounded variation and no drift on a compact interval $I$ of $\mathbb{R}$. It is a standard result that with probability $1$, $X^*$ reaches its maximum at only one point of $I$, that we denote by $T$.

\begin{Lm}
\label{Dissym2}
If $T$ is a positive jump time, then $0$ is irregular for the positive half-line, and $T$ is isolated on the left but not on the right.\end{Lm}

\begin{proof}

Let us assume that $T$ is a positive jump time. It is  a stopping time because it is the time of crossing of  $z \in \mathbb{Q} \cap ]\sup_{s<T}X(s),X({T})[$. Since $X^*(T)$ is only reached once, (in $T$,) $\mathbb{P}(R_{T}^+(X))=0=p^+(X)$, whence $0$ is irregular for the positive half-line. Moreover, $\lim \limits_{a \to 0^+}\frac{X(T+a)-X(T)}{a} =0$ because $T$ is a stopping time, whence $T$ is not right isolated. $T$ is left isolated as a positive jump.
\end{proof}
We have the converse:
\begin{Lm}
If $0$ is irregular for the positive half-line, then a.s. $T$ is a positive jump time.
\end{Lm}

\begin{proof}
Let $V_{1}=\inf\{a >0~;~X(a)\geq 0\}$ and $V_{n+1}=V_{1}\circ \theta_{V_{n}}+V_{n}$, it is a sequence of stopping times with iid increments.  If $V_{1}=+\infty$, simply set $V_{n}=+\infty$ for $n\geq 1$. The point $0$ being irregular for the positive half-line,  the times $V_{n}$ are a.s. strictly positive.
In particular, a.s. $\lim_{n}V_{n}=\infty$ and the times $V_{n}$ are discrete in $\mathbb{R}^+$.

Let us also show that $X_{V_{1}}>0$.
It suffices to notice that if $X_{V_{1}}=0$, if we restrict ourselves to an interval of the form $[0,q]$, $q\in \mathbb{Q}\cap ]V_{1},V_{2}[,$ $X$ reaches twice its maximum on this interval, which almost never happens.
$V_{1}$ is then a positive jump time, and so are the $V_{n}$.

It is clear that $T=V_{n+1}$, where $n$ is the greatest integer such that $X(V_{n})<T$, and $T$ is therefore a positive jump time.
\end{proof}
We can conclude the demonstration of the theorem by this last lemma:
\begin{Lm}
\label{PasDissym}
We suppose that  $0$ is regular for the positive and negative half-lines.
Then  $\overline{X}'(T^-)=\overline{X}'(T)=0$ a.s., and $T$ is isolated neither on its right nor on its left.
\end{Lm}

\begin{proof}
Let $\Omega_{1}$ be the event $\overline{X}'(T^-)>0$.
We show by contradiction that $\Omega_{1}$ is negligible.

Suppose without loss of generality that $I=[0,1]$ and that $X$ is defined on all the real line. Denote by $X_{[a,b]}$ the restriction of $X$ on a segment $[a,b]$, and $T_{[a,b]}$ the time where $X^*_{[a,b]}$ reaches its maximum. Consider the following event
\begin{equation}
\Omega_{2}=(\Omega_{1};~T_{[0,1]}<T_{[0,2]};~T_{[0,1]}\in \Ext^+_{X_{[0,2]}}).
\end{equation} 

The event $\Omega_{2}$ is realised if $T_{[0,1]}$ is an extremal superior time for $X_{[0,2]}$, while being strictly inferior to $ T_{[0,2]}$ where the process reaches its maximum. Hence according to Proposition \ref{PpsJumpPext}, $T_{[0,1]}$ is a positive jump time for $X_{[0,2]}$, thus for $X_{[0,1]}$ as well. This almost never happens according to \ref{Dissym2} since we assumed that $\{0\}$ was regular for the positive half-line. We will hence arrive at a contradiction if we show that $\Omega_{2}$ has positive probability.

Assume that for some $\beta>0, m>0$ and $s\in [0,1)$ we have
\begin{equation}
\label{eq:area}
\begin{cases}
T_{[0,1]}\leq s,\\
\overline{X_{[0,1]}}'(T_{1}^-)>\beta,\\
m<X(T_{[0,1]})<m+(1/2)\beta(1-s),\\
m+(1/2)\beta(1-s)<X(T_{[0,2]})<m+\beta(T_{[0,2]}-s),
\end{cases}
\end{equation}
then $T_{[0,1]}$ is in $\Ext^+_{X_{[0,2]}}$. Indeed, the slope of the line going through $M_{T_{[0,1]}}$ and $M_{T_{[0,2]}}$ satisfies
\begin{equation*}
\frac{X(T_{[0,2]})-X(T_{[0,1]})}{T_{[0,2]}-T_{[0,1]}}<\frac{m+\beta(T_{[0,2]}-s)-X(T_{[0,1]})}{T_{[0,2]}-s}<\beta<\overline{X_{[0,1]}}'(T_{[0,1]}^-).
\end{equation*}

All the problem amounts to finding $\beta,s,m$ such that the event (\ref{eq:area}) has positive probability. Let us fix $\beta>0$ such that the event $\Omega_{\beta}=(X(T_{1}^-)>\beta)$ has positive probability, and let $s<1$ be such that $\mathbb{P}(\Omega_{\beta};T_{[0,1]}\leq s)>0$. Conditionally 
on this event, there is $m>0$ such that $$\Omega_{\beta,m,s}=(\Omega_{\beta};  T_{[0,1]}\leq s;m<X(T[0,1])<m+1/2 \beta(1-s); X(1)>m)$$ has positive probability.

 Since $X$ is a Markov process, the trajectory of $X_{[0,2]}$ is independent of that of $X_{[0,1]}$ given $X(1)$. Thus, conditionally on $\Omega_{\beta,m,s}$, $X_{[0,2]}$'s graph (which left end-point is $(1,X(1))$ is comprised in the rectangle $[1,2]\times (m,m+\beta(1-s))$ with positive probability. Also, the maximal value of $X_{[0,2]}$ (i.e. $X(T_{[0,2]})$) is larger than $m+(1/2)\beta(1-s)$ with positive probability. Thus, with these particular values, (\ref{eq:area}) occurs with positive probability.

We hence proved by contradiction  $\overline{X}'(T^-)=0$ a.s., and by reflexivity we also have $\overline{X}'(T)=0$. Since $T$ is the only point where $X$ reaches its maximum, this implies that $T$ is a right and left accumulation point of $E_{X}^+$.
\end{proof}

\section*{Acknowledgements}
I wish to thank Youri Davydov, my PhD advisor, at the origin of this paper, for his guidance in this work. The comments and remarks of both referees have also been very important in the improvement of the paper.

\end{document}